\documentclass[11pt]{amsart}
\usepackage{geometry}                
\usepackage{graphicx}
\usepackage{amssymb}
\usepackage{epstopdf}
\usepackage{caption}
\usepackage{subcaption}
\usepackage{enumerate}
\newcommand{\lline}{\overleftrightarrow}
\newcommand{\ppoint}{\widehat}
\usepackage{url}
\newtheorem{thm}{Theorem}

\title{A Six-Point Ceva-Menelaus Theorem}
\author[B.D.S. ``Blue'' McConnell]{B.D.S. ``Blue'' McConnell \\ \tt{math@daylateanddollarshort.com}}
\date{\today}

\begin{document}
\begin{abstract}
We provide a companion to the recent B{\'e}nyi-{\'C}urgus generalization of the well-known theorems of Ceva and Menelaus, so as to characterize both the collinearity of points and the concurrence of lines determined by six points on the edges of a triangle. A companion for the generalized area formula of Routh appears, as well.
\end{abstract}
\maketitle
The venerable theorems of (Giovanni) Ceva and Menelaus (of Alexandria) concern points on the edge-lines of a triangle. Each point defines ---and is defined by--- the ratio of lengths\footnote{Throughout, we consider segment lengths to be {\em signed}, with each of $\overrightarrow{AB}$, $\overrightarrow{BC}$, $\overrightarrow{CA}$ ---for {\em distinct} $A$, $B$, $C$--- indicating the direction of a positively-signed segment on the corresponding (extended) side of the triangle. Moreover, we adopt these conventions regarding ratios of these lengths:
\begin{equation*}
\frac{|PP|}{|PQ|} = 0 \qquad\qquad \frac{|PQ|}{|QQ|} = \infty \qquad\qquad \frac{|PX|}{|XQ|} = -1, \;\; \text{for $X$ the point at infinity on $\lline{PQ}$}
\end{equation*}} of collinear segments joining it to two of the triangle's vertices, and the theorems use a trio of such ratios to neatly characterize the special configurations in Figure 1.

\vspace{-0.7em}
\begin{figure}[h]
\centering
\begin{subfigure}[t]{0.5\textwidth}
\centering
\includegraphics[width=.8\textwidth]{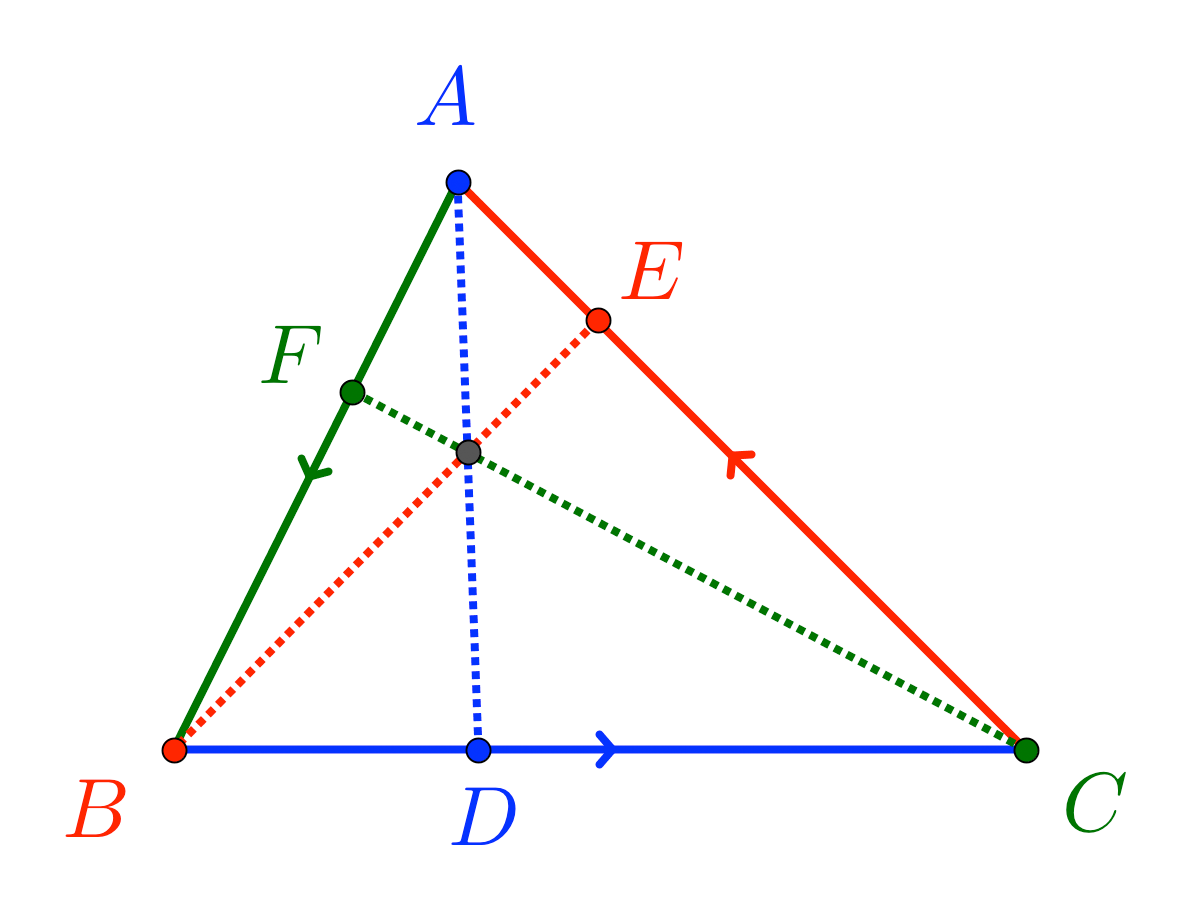}
\caption[justification=centerlast]{Lines through a triangle's vertices, \\ meeting at a common point. (Ceva)}
\end{subfigure}\hfill
\begin{subfigure}[t]{0.5\textwidth}
\centering
\includegraphics[width=.8\textwidth]{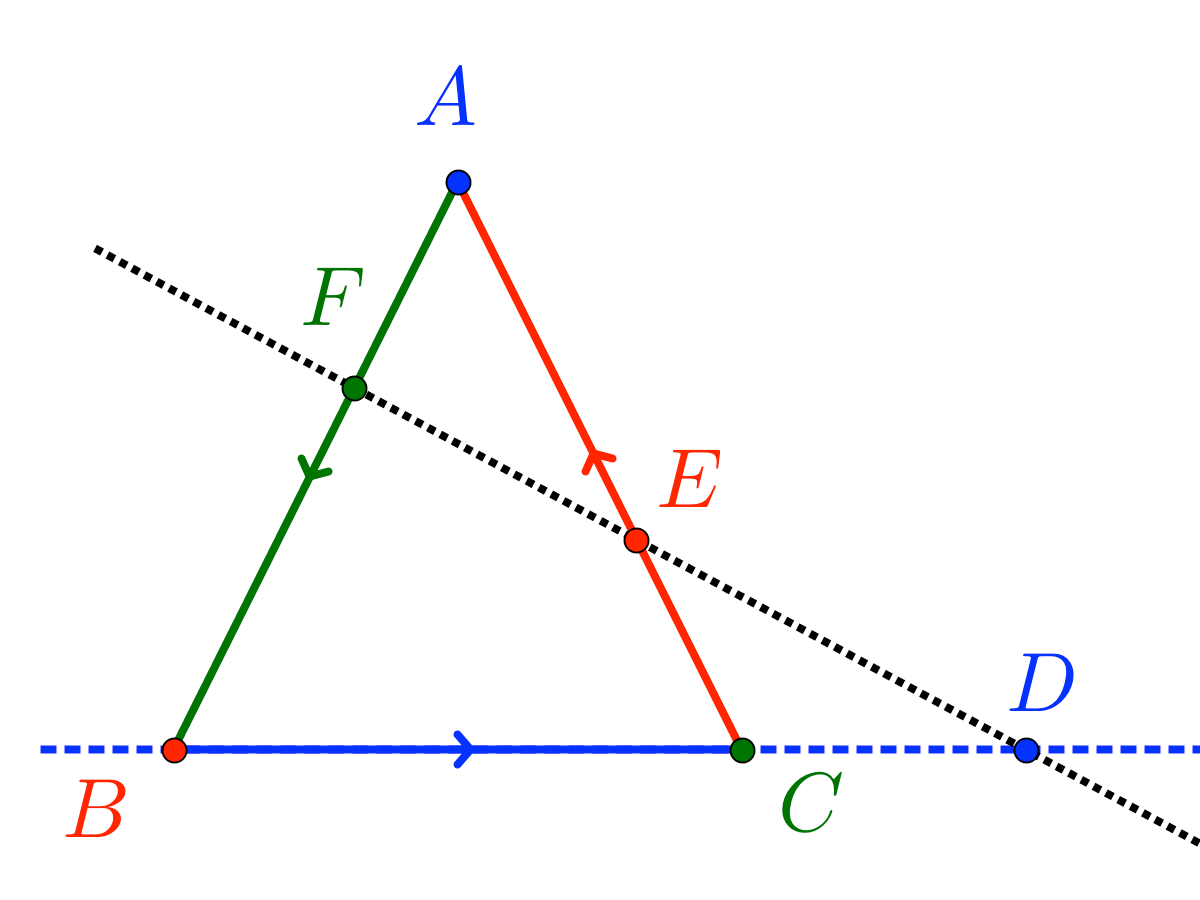}
\caption[justification=centerlast]{Points on a triangle's edges, \\ lying on a common line. (Menelaus)}
\end{subfigure}
\caption{}
\end{figure}

\vspace{-0.2em}
Specifically, with $D$, $E$, $F$ on edge-lines opposite respective vertices $A$, $B$, $C$, we write
\begin{equation}
d := \frac{|BD|}{|DC|} \qquad e := \frac{|CE|}{|EA|} \qquad f := \frac{|AF|}{|FB|}
\end{equation}
and express the theorems as follows:
\begin{subequations}
\begin{thm}[Ceva]
Lines $\overleftrightarrow{AD}$, $\overleftrightarrow{BE}$, $\overleftrightarrow{CF}$ pass through a common point if and only if \begin{equation}
\label{eqn_ceva}
d e f = \phantom{-}1
\end{equation}
\end{thm}
\begin{thm}[Menelaus]
Points $D$, $E$, $F$ lie on a common line if and only if
\begin{equation}
\label{eqn_menelaus}
d e f = -1
\end{equation}
\end{thm}
\end{subequations}

\vspace{1em}
B{\'e}nyi and {\'C}urgus \cite{BC2012}, and this author, independently (and nearly-simultaneously) considered separate aspects of the same approach to generalizing the above ---namely, doubling the number of points on the triangle's edges--- arriving at equations whose terms, in the grand Ceva-Menelaus tradition, differ only in sign.

Interestingly, the B{\' e}nyi-{\' C}urgus result concerns {\em Ceva}-like elements (lines through vertices) and a {\em Menelaus}-like phenomenon (collinearity of points). This author's contribution, on the other hand, concerns {\em Menelaus}-like elements (points on edges) and a {\em Ceva}-like phenomenon (concurrence of lines). 

\begin{figure}[h]
\centering
\begin{subfigure}[t]{0.5\textwidth}
\centering
\includegraphics[width=.8\textwidth]{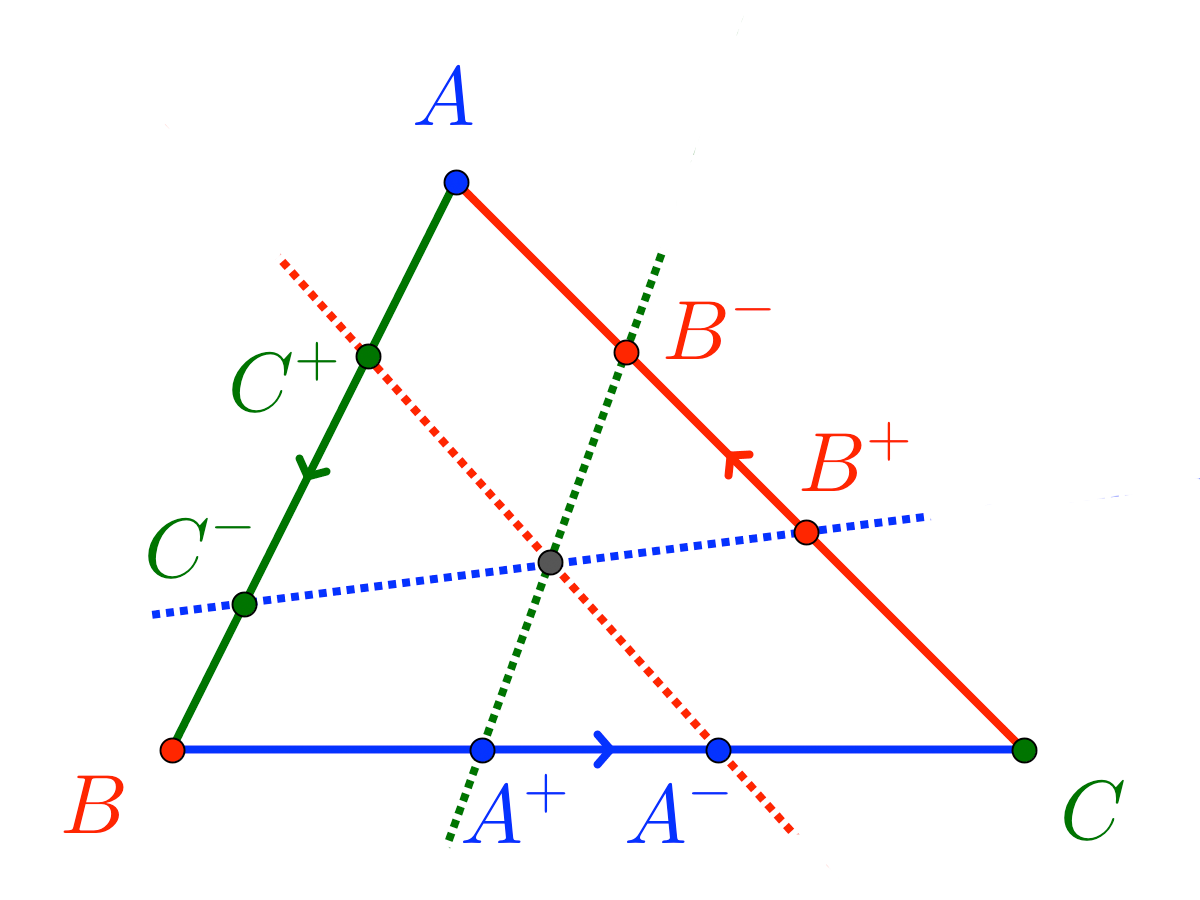}
\caption[justification=centerlast]{Lines through points on edges, \\ meeting a common point.}
\end{subfigure}\hfill
\begin{subfigure}[t]{0.5\textwidth}
\centering
\includegraphics[width=.8\textwidth]{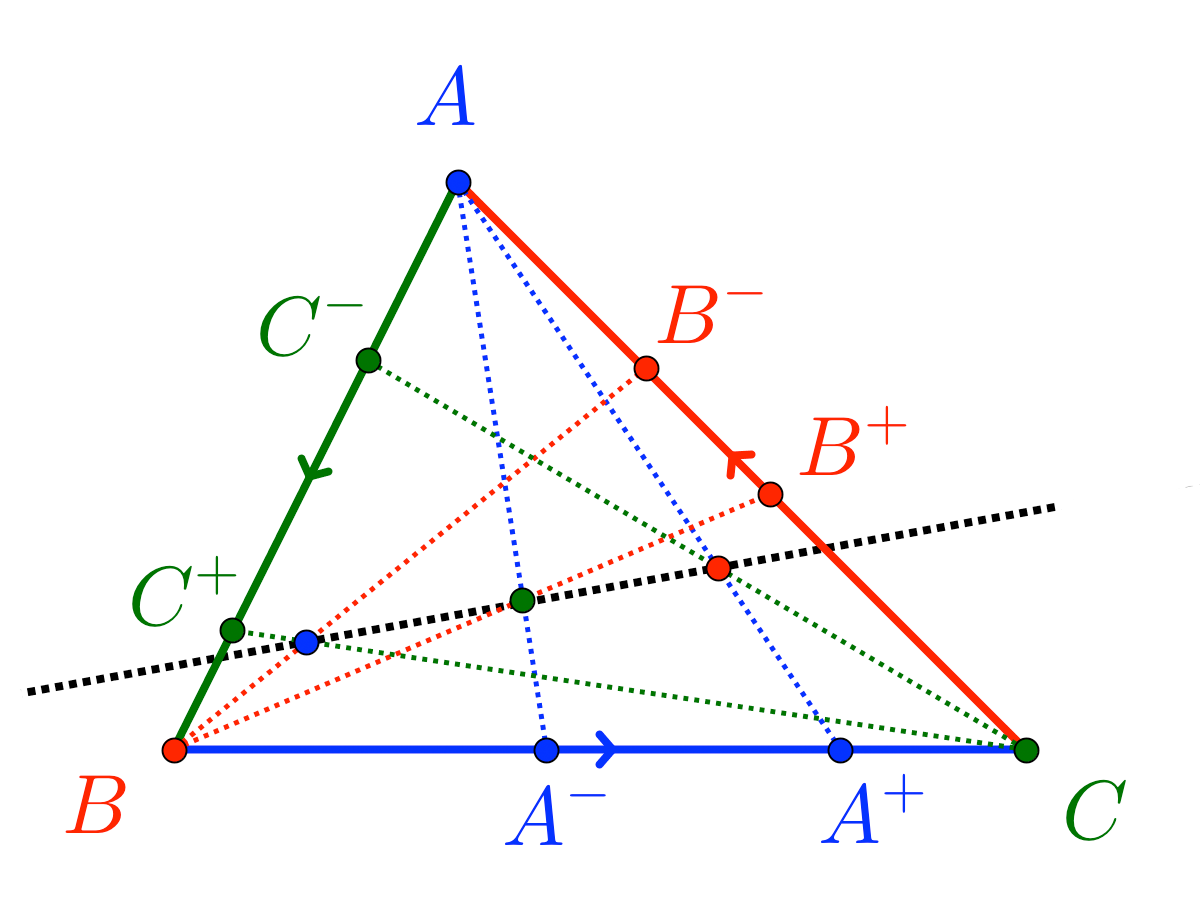}
\caption[justification=centerlast]{Points on lines through vertices, \\ lying on a common line. (B{\'e}nyi-{\'C}urgus)}
\end{subfigure}
\caption{}
\end{figure}

Place points $A^{+}$ and $A^{-}$ on the edge-line opposite vertex $A$; likewise, $B^{+}$ and $B^{-}$ opposite $B$, and $C^{+}$ and $C^{-}$ opposite $C$. Define these ratios:\footnote{Observe that the superscripts emphasize an opposing directionality in the definitions of the ratios. For instance, the points in ratio $a^{+}$ trace the path $B$-$A^{+}$-$C$, with endpoints oriented in the {\em positive} direction; in $a^{-}$, the path $C$-$A^{-}$-$B$ has endpoints oriented in the {\em negative} direction. Were we to define all six ratios in ``matching'' orientations ---as was done in \cite{BC2012}--- the resulting formulas would lose some clarity and symmetry.} 

\begin{equation}\begin{array}{c}
\displaystyle a^{+} := \frac{|BA^{+}|}{|A^{+}C|} \qquad b^{+} := \frac{|CB^{+}|}{|B^{+}A|} \qquad c^{+} := \frac{|AC^{+}|}{|C^{+}B|} \\[10pt]
\displaystyle a^{-} := \frac{|CA^{-}|}{|A^{-}B|} \qquad b^{-} := \frac{|AB^{-}|}{|B^{-}C|} \qquad c^{-} := \frac{|BC^{-}|}{|C^{-}A|}
\end{array}
\end{equation}

\pagebreak
\begin{thm}[Six-Point Ceva-Menelaus Theorem] $\phantom{xyzzy}$
\begin{subequations}
\begin{enumerate}[\hspace{1em}a.]
\item Lines $\lline{B^{+}C^{-}}$, $\lline{C^{+}A^{-}}$, $\lline{A^{+}B^{-}}$ pass through a common point if and only if
\begin{equation}
\label{eqn_6ptconcurrence}
a^{+} b^{+} c^{+} \;+\; a^{-} b^{-} c^{-} \;=\; \phantom{-}1 \;-\; a^{+} a^{-} \;-\; b^{+} b^{-} \;-\; c^{+} c^{-}
\end{equation}
\item (B{\'e}nyi-{\'C}urgus) Points\footnote{To amplify the duality with (a), we write ``$\ppoint{P^{-}Q^{+}}$'' for the point of intersection of lines $\lline{PP^{-}}$ and $\lline{QQ^{+}}$.} $\ppoint{B^{-}C^{+}}$, $\ppoint{C^{-}A^{+}}$, $\ppoint{A^{-}B^{+}}$ lie on a common line if and only if
\begin{equation}
\label{eqn_6ptcollinearity}
a^{+} b^{+} c^{+} \;+\; a^{-}b^{-}c^{-} \;=\; -1 \;+\; a^{+} a^{-} \;+\; b^{+} b^{-} \;+\; c^{+} c^{-}
\end{equation}
\end{enumerate}
\end{subequations}
\end{thm}

Note: Identifying $A^{-}$, $B^{-}$, $C^{-}$ with $C$, $B$, $A$ yields $a^{-} = b^{-} = c^{-} = 0$, so that (\ref{eqn_6ptconcurrence}) and (\ref{eqn_6ptcollinearity}) reduce to (\ref{eqn_ceva}) and (\ref{eqn_menelaus}). The Six-Point Theorem generalizes the traditional results.

For proof, one can invoke vector techniques, as indicated with Theorem \ref{thm_sixpointrouth} below.

\subsection*{Routh, too.} When Ceva's lines fail to concur, and when Menelaus' points fail to ``colline'', they determine triangles. One might well ask how the area\footnote{As with length, we consider triangle area to be {\em signed}. Areas $|\triangle ABC|$ and $|\triangle DEF|$ agree in sign when vertex paths $A$-$B$-$C$-$A$ and $D$-$E$-$F$-$D$ trace their respective figures in the same direction; similarly for triangles defined by their edge-lines.} of each resulting triangle compares to that of the original figure. (Edward John) Routh provided answers.  (See \cite{BC2012}.)

\vspace{-0.5em}
\begin{figure}[h]
\begin{subfigure}[t]{0.5\textwidth}
\centering
\includegraphics[width=.8\textwidth]{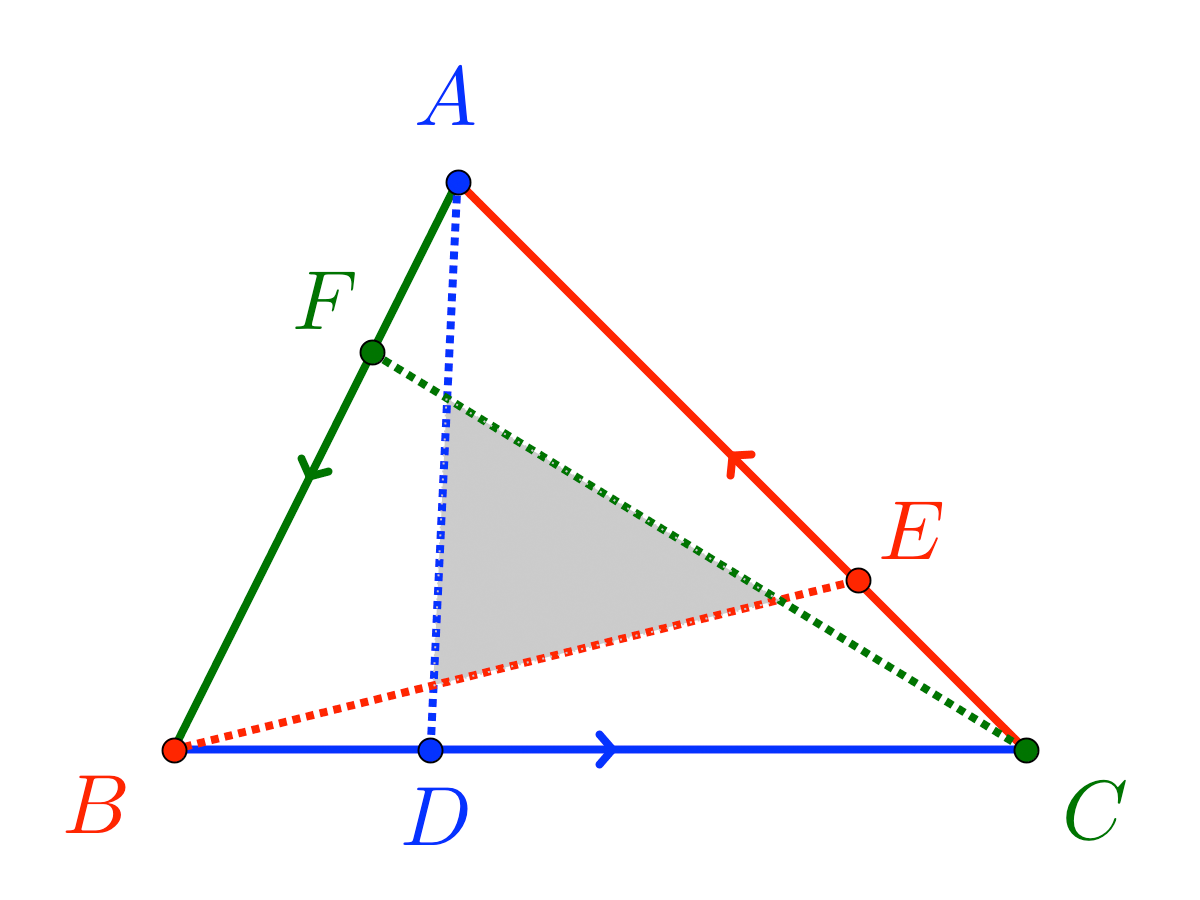}
\caption[justification=centerlast]{Triangle from lines through vertices.}
\end{subfigure}\hfill
\begin{subfigure}[t]{0.5\textwidth}
\centering
\includegraphics[width=.8\textwidth]{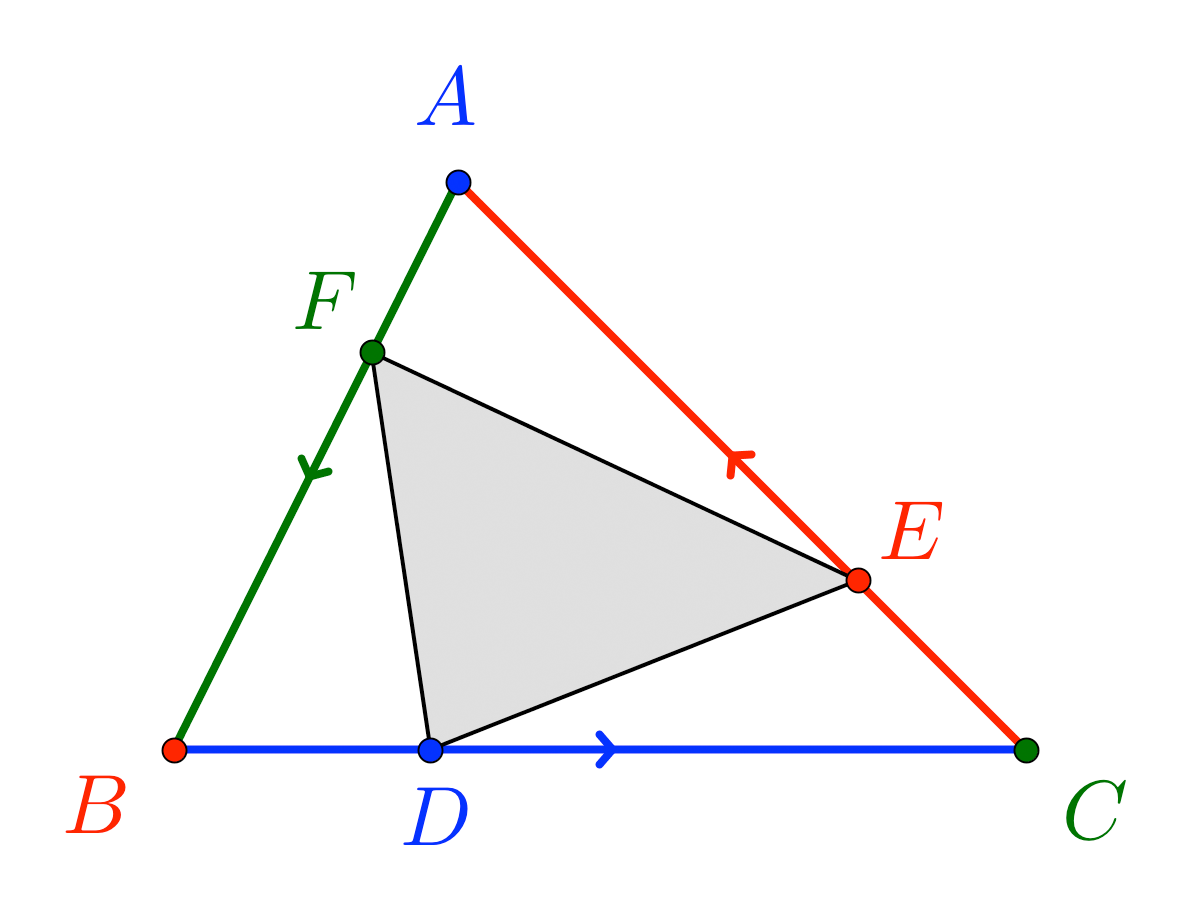}
\caption[justification=centerlast]{Triangle from points on edges.}
\end{subfigure}
\caption{}
\end{figure}

\begin{thm}[Routh's Formulas] \hspace{3em}
\begin{subequations}
\begin{enumerate}[\hspace{1em}a.]
\item (``Routh's Theorem''). The triangle with (non-parallel) edge-lines $\lline{AD}$, $\lline{BE}$, $\lline{CF}$ has area
\begin{equation}
\label{eqn_routh2}
|\triangle ABC| \cdot \frac{\left(\; d e f - 1 \;\right)^2}{
    \left(\; 1 + d + d e \;\right)\left(\; 1 + e + e f \;\right) \left(\; 1 + f + f d \;\right)}
\end{equation}
\item The triangle with (finite) vertices $D$, $E$, $F$ has area \\
\begin{align}
\label{eqn_routh1}
|\triangle ABC| \cdot \frac{\; d e f + 1 \;}{
    \left(\; 1 + d \;\right) \left(\; 1 + e \;\right) \left(\; 1 + f \;\right)}
\end{align}
\end{enumerate}
\end{subequations}
\end{thm}

Observe that the numerator in each of these formulas ---and, thus, the area of the triangle in question--- vanishes, as it should, when (and only when) the conditions for Ceva's or Menelaus' theorems indicate that the triangle degenerates into a point or a line. The reader may verify that a denominator vanishes when (and only when) the triangle becomes unbounded, having non-finite vertices and parallel edges.

\vspace{1em}
B{\'e}nyi and {\'C}urgus \cite{BC2012} specifically address the six-point generalization of Theorem 4b. At the suggestion of Mr. {\'C}urgus, this author derived the counterpart generalization of 4a.

\vspace{-0.5em}
\begin{figure}[h]
\begin{subfigure}[t]{0.5\textwidth}
\centering
\includegraphics[width=.8\textwidth]{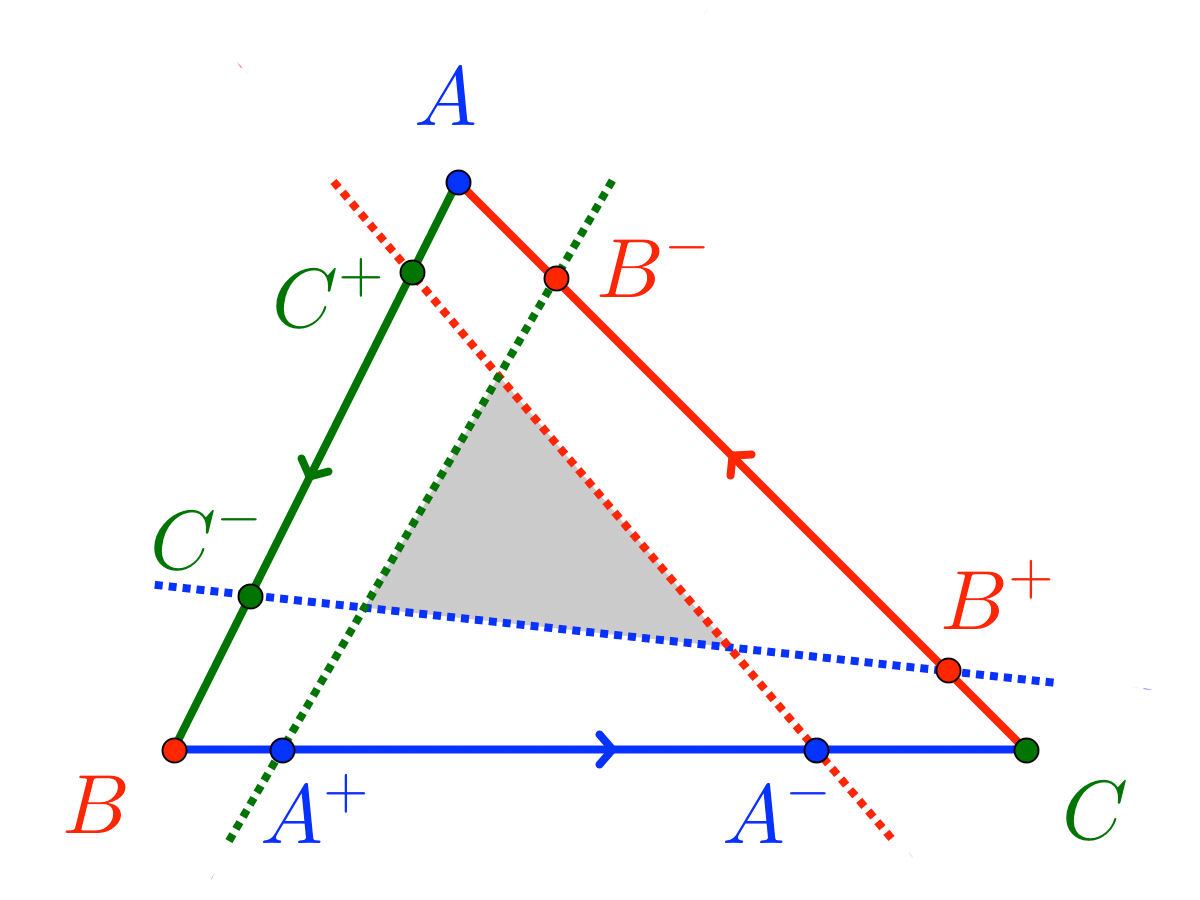}
\caption[justification=centerlast]{Triangle from lines through points on edges.}
\end{subfigure}\hfill
\begin{subfigure}[t]{0.5\textwidth}
\centering
\includegraphics[width=.8\textwidth]{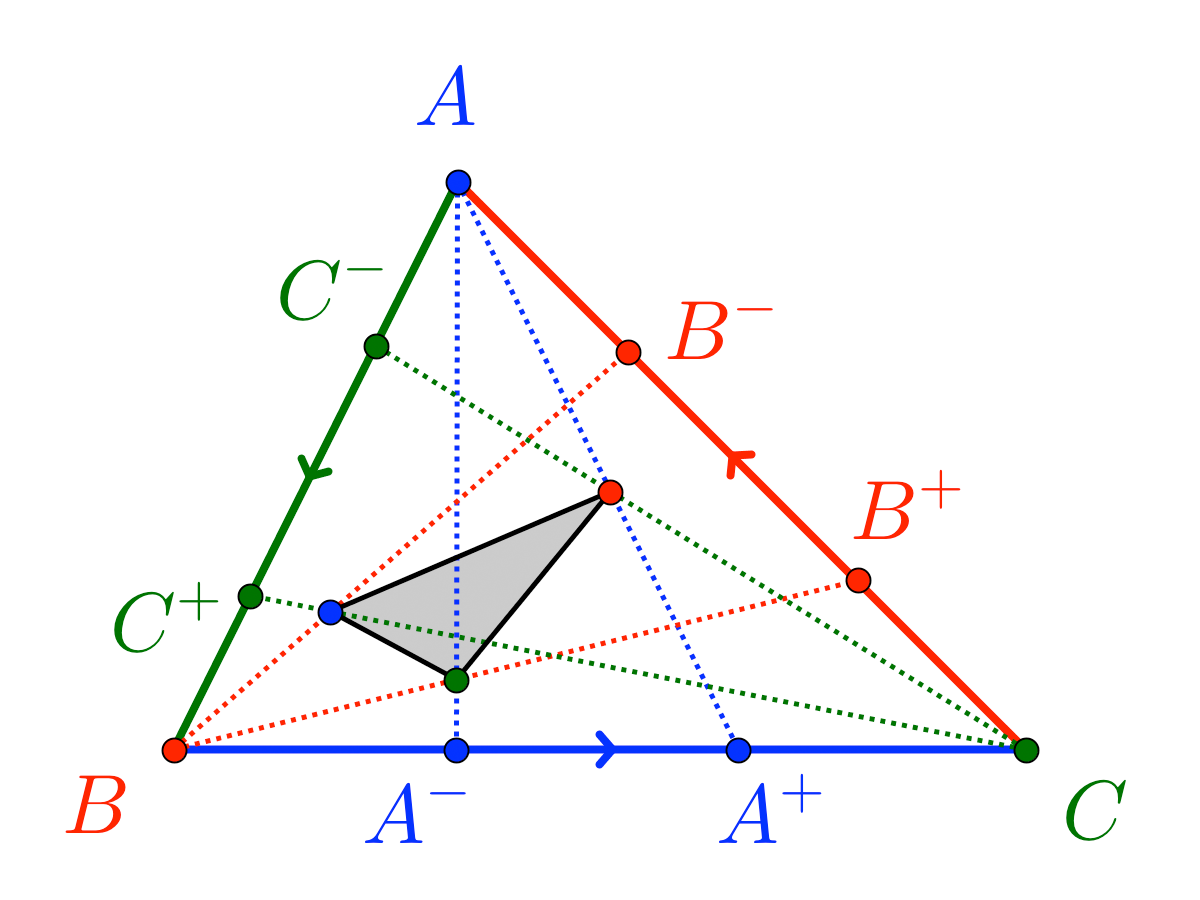}
\caption[justification=centerlast]{Triangle from points on lines through vertices. (B{\'e}nyi-{\'C}urgus)}
\end{subfigure}
\caption{}
\end{figure}

\begin{thm}[Six-Point Routh Formulas] $\phantom{xyzzy}$
\label{thm_sixpointrouth}
\begin{subequations}
\begin{enumerate}[\hspace{1em}a.]
\item The triangle with (non-parallel) edge-lines $\lline{B^{+}C^{-}}$, $\lline{C^{+}A^{-}}$, $\lline{A^{+}B^{-}}$ has area \\
\begin{align}
\label{eqn_6ptrouth1}
|\triangle ABC| \cdot \frac{\left(\; a^{+} b^{+} c^{+} + a^{-} b^{-} c^{-} + a^{+} a^{-} + b^{+} b^{-} + 
    c^{+} c^{-} - 1  \;\right)^2}{
    \begin{array}{c}
    \phantom{\cdot}\left(\; 1 - a^{+} a^{-} + b^{-} ( 1 + a^{-} ) + c^{+} ( 1 + a^{+} ) \;\right) \\[3pt]
    \cdot\left(\; 1 - b^{+} b^{-} + c^{-} (  1 + b^{-} ) + a^{+} ( 1 + b^{+} )  \;\right) \\[3pt]
    \cdot\left(\; 1 - c^{+} c^{-} + a^{-} ( 1 + c^{-} ) + b^{+} ( 1 + c^{+} )  \;\right)
    \end{array}}
\end{align}
\item (B{\'e}nyi-{\'C}urgus). The triangle with (finite) vertices $\ppoint{B^{-}C^{+}}$, $\ppoint{C^{-}A^{+}}$, $\ppoint{A^{-}B^{+}}$ has area
\begin{equation}
\label{eqn_6ptrouth2}
|\triangle ABC| \cdot \frac{\; a^{+} b^{+} c^{+} + a^{-} b^{-} c^{-} - a^{+} a^{-} - b^{+} b^{-} - 
    c^{+} c^{-} + 1 \;}{
    \left(\; 1 + b^{-} + c^{+} \;\right)\left(\; 1 + c^{-} + a^{+} \;\right) \left(\; 1 + a^{-} + b^{+} \;\right)}
\end{equation}
\end{enumerate}
\end{subequations}
\end{thm}

\vspace{1em}
\begin{proof} Treating points as vectors, we can write
\begin{equation}
A^{+} = \frac{B (1+a^{+}) + Ca^{+}}{1+a^{+}} \qquad A^{-} = \frac{Ba^{-} + C(1+a^{-})}{1+a^{-}} \qquad \text{, etc.}
\end{equation}
to find, after a bit of tedious algebra, that the vertices of the triangles in parts (a) and (b) of the theorem have the respective forms
\begin{equation}
\frac{A ( 1 - a^{-} a^{+}) + B ( c^{+} + a^{-} b^{-} ) + C ( b^{-} + a^{+} c^{+} )}
{1 - a^{-} a^{+} + b^{-} ( 1 + a^{-} ) + c^{+} ( 1 + a^{+} )}
\qquad\text{and}\qquad
\frac{A + B c^{+} + Cb^{-}}{1 + c^{+} + b^{-}}
\end{equation}
The area formulas follow from a bit more ---and more-tedious--- algebra. Of course, since ratios of lengths of collinear segments, and of areas of coplanar triangles, are preserved under affine transformation, one could simplify this analysis somewhat by assuming, say, $A = (0,0)$, $B=(1,0)$, $C=(0,1)$; even in generality, however, verification of these formulas amounts to just a few seconds' effort from a computer algebra system.
\end{proof}

\vspace{2em}
\subsection*{Remarks} We can accentuate the duality of the traditional results of Ceva and Menelaus by reciting them thusly: ``{\em Points determined by pairs of lines through the (vertex-)points of a triangle coincide if and only if ...}'' versus ``{\em Lines determined by pairs of points on the (edge-)lines of a triangle coincide if and only if ...}''.

Taking a deep breath, we can do likewise for the parts of the Six-Point Theorem: ``{\em Points determined by pairs of lines determined by pairs of points on the (edge-)lines of a triangle coincide if and only if ...}'' versus  ``{\em Lines determined by pairs of points determined by pairs of lines through the (vertex-)points of a triangle coincide if and only if ...}''.

What formulas characterize the coincidence of lines and/or points at the next order of complexity? For that matter, what strategy for pairing lines and/or points best {\em constitutes} the next order of complexity?

\end{document}